\documentclass[12pt,letterpaper, final]{amsart}
\usepackage{amssymb, amsmath, amsthm}
\usepackage{mathrsfs}  
\usepackage[dvipsnames]{xcolor}
\usepackage{graphicx}
\usepackage[colorlinks=true, citecolor=blue, linkcolor=WildStrawberry, urlcolor=blue]{hyperref}
\usepackage[text={6in,8in},centering]{geometry}
 \usepackage{xcolor}
\usepackage{mathrsfs}  
\usepackage{bm}
\usepackage{ifdraft}
\ifoptionfinal{
\usepackage[disable]{todonotes}
}{
\usepackage[norefs, nocites]{refcheck}

\usepackage[notref, notcite]{showkeys}
\usepackage[bordercolor=white, color=white]{todonotes}
}

\makeatletter\providecommand\@dotsep{5}\def\listtodoname{List of Todos}\def\listoftodos{\hypersetup{linkcolor=black}\@starttoc{tdo}\listtodoname\hypersetup{linkcolor=blue}}\makeatother
 
\usepackage{hyperref}
\usepackage{etoolbox}
\newtheorem{theorem}{Theorem}[section]

\theoremstyle{definition}

\theoremstyle{remark}
\newtheorem{remark}[theorem]{Remark}

\numberwithin{equation}{section}

\renewcommand{\leq}{\leqslant}
\renewcommand{\geq}{\geqslant}

\newcommand*\xbar[1]{%
   \hbox{%
     \vbox{%
       \hrule height 0.5pt 
       \kern0.5ex
       \hbox{%
         \ensuremath{#1}%
       }%
     }%
   }%
}

\makeatletter
\def\@keywordsname{Key words}
\makeatother

\title[]{\textbf{Anisotropic Calder\'{o}n problem of a nearly Laplace-Beltrami operator of order $2+$}}
\author[ ]{\small  Susovan Pramanik }
\date{}

  \address{Susovan Pramanik, Harish-Chandra Research Institute, A CI of Homi Bhabha National Institute, Chhatnag Road, Jhunsi, Allahabad 211 019, India}
\email{susovanpramanik@hri.res.in}

\keywords{Inverse problems, log-Laplacian, Anisotropic Calder\'on problem, Symbol class of  operators, Functional Calculus}

\begin{document}

\maketitle
\noindent
\begin{abstract}
This paper investigates the anisotropic Calder\'{o}n problem for Logarithemic Laplacian, on closed Riemannian manifolds, which could be considered as near  Laplace operator. We demonstrate that the Cauchy data set recovers the geometry of a closed Riemannian manifold up to standard gauge.  
\end{abstract}

\section{Introduction}

The anisotropic Calder\'{o}n problem is geometric in nature (cf. \cite{LU89, Sal13}). Let us state the well-known open problem in the context of a closed (compact without boundary), connected Riemannian manifold $(M, g)$ with $\dim(M) > 2$. Let $\mathcal{O} \subset M$ be a non-empty open set. The Calder\'{o}n problem examines whether knowledge of the Cauchy data set 
\begin{equation}\label{0pen}
    \widetilde{\mathcal{C}}^\mathcal{O}_{M, g}= \{(u|_\mathcal{O}, (-\Delta_g)u|_\mathcal{O}) \mid u \in C^\infty(M), \, -\Delta_g u = 0 \,\text{in} \,\, M\setminus \overline{\mathcal{O}}\}
\end{equation}
determines the isometry class of the manifold $(M, g)$. \((-\Delta_g) \) represents the positive Laplace-Beltrami operator on \( (M, g)\). 

\medskip
When $\dim(M)= 2$, the problem is solved using an additional gauge due to conformal invariance of the Laplace-Beltrami operator. See some of the references \cite{Nach96, KL06, Buk08, GT11}. The problem for real-analytic manifolds in dimensions three and above has been addressed (see the references \cite{LU89,  LU01, LTU03}), but it is still open for smooth manifolds. Positive results have been observed in transversally anisotropic geometries, see \cite{DKSU09, DKLS16}. Starting from these celebrated works \cite{SU87, KSU07}, for current research on the Calder\'{o}n problem, we refer to the survey article \cite{Uhl14}.

\medskip
Banking on the fact that, for all $s>0$, $\lim_{\lambda\to\infty}\frac{\ln\lambda}{\lambda^s}=0$, as an unbounded operator $-\Delta_g \circ \ln[(-\Delta_g)]$, which we would call a logarithmic Laplacian operator, can be considered as a nearly Laplace-Beltrami operator $(-\Delta_g)$ of order $2+$. We are interested in addressing the corresponding anisotropic Calder\'{o}n problem, which is related to the non-local nature of the logarithmic Laplacian operator. 

\medskip
To be more specific, let us set some $m> 1$ and consider the following operator on $(M,g)$ as:
\begin{equation}\label{log-Laplace} 
\mathcal{L}_g:=(-\Delta_g+m\mathbb{I}) \circ   \ln [(-\Delta_g+m\mathbb{I})]= \ln [(-\Delta_g+m\mathbb{I})]\circ (-\Delta_g+m\mathbb{I}).
\end{equation}
The purpose of perturbing the Laplace operator by $m\mathbb{I}$ where $m>1$ to ensure the operator $\mathcal{L}_g$ remains elliptic.

\medskip
The goal of this work is to explore the anisotropic Calder\'{o}n problem, whether the following Cauchy data set 
\[ \mathcal{C}^{\mathcal{O}}_{M, g}= \{(u|_\mathcal{O},  \mathcal{L}_g u|_\mathcal{O}) \mid u \in C^\infty(M), \,\,\, \mathcal{L}_g u = 0 \,\text{in} \, M \setminus \overline{\mathcal{O}} \} \]  
determines the metric $g$ modulo the standard gauge. This is the inverse problem we are interested in. 

\medskip
In a different set up, recently the authors here \cite{HLW24}, looked into inverse problems of recovering the potential, associated to the near-zero logarithmic operator modelled into Schr\"{o}dinger-type equations in flat spaces.  

\medskip
The motivation of our work is driven by the works see \cite{FGKU25, FGKRSU25} for both closed and open manifolds, respectively, on the nonlocal Calder\'{o}n's problem of recovering the anisotropic medium or the metric, where the nonlocal operator is modeled by the fractional Laplace-Beltrami operator of fixed order $2\alpha$, $\alpha\in (0,1)$. We now know that the Cauchy data set 
\begin{align*}
    \mathcal{C}^{\alpha, \mathcal{O}}_{M, g}= \{(u|_\mathcal{O}, (-\Delta_g)^\alpha u|_\mathcal{O}) \mid u \in C^\infty(M), \, (-\Delta_g)^\alpha u = 0 \,\text{in} \,\, M\setminus \overline{\mathcal{O}}, \,\, \alpha\in (0,1)\}
\end{align*}
uniquely specifies the manifold $(M, g)$ up to a diffeomorphism. The results stem from the work \cite{GU21}, and subsequent development adds \cite{CGRU23, Feiz24, CO24, FGKU25, FKU24, FGKRSU25, CR25} to further complete the picture in this nonlocal regime. 

\medskip
 This is based on reducing the problem into the hyperbolic inverse problem, and thanks to the work of \cite{HLOS18} for instance, we have a definite conclusion in the nonlocal analogue of the anisotropic Calder\'{o}n problem, which is still to be discovered for the local anisotropic Calder\'{o}n problem.

\medskip
Recently, we orchestrated such exercise in solving the inverse problem of recovering the geometry of closed manifold $(M, g)$ from the Cauchy data associated with the non-local operator $((-\Delta_g)^2+m^2\mathbb{I})^{\frac{1}{2}}$ of order exactly $2$, see \cite{SP25}.

\medskip
In this article, our non-local operator is modelled by the logarithmic Laplacian operator (cf. \eqref{log-Laplace}), which is of order $2+$ and dynamically closer to the open problem linked to the Laplace operator of order $2$ (cf. \eqref{0pen}), rather than the fractional Laplace operator of some fixed order less than $2$. Our work illustrates that the method of proof indicated in the preceding paragraph can be extended to include our nonlocal operator, so as to solve the inverse problem globally.

\subsection*{Formal Definitions and Main Theorem}
Let $(M, g)$ be a closed, connected Riemannian manifold with dimension \(\dim(M) \geq 2 \). We define \(-\Delta_g \) as the positive definite, self-adjoint Laplace-Beltrami operator on \(L^2(M) \) with domain $\mathcal{D}(-\Delta_g) = H^2(M);$ \cite[pp. 25]{Tay23}.

The spectrum of \( (-\Delta_g) \) over $(M, g)$ consists of a discrete sequence of eigenvalues  
\[
0 = \lambda_0 < \lambda_1 \leq \lambda_2 \leq \cdots \to+\infty,
\]  
where each eigenvalue \( \lambda_k \) has a finite multiplicity \( d_k \).

\medskip
Define \(\pi_k: L^2(M) \to \ker(-\Delta_g - \lambda_k) \) as the orthogonal projection onto the eigenspace corresponding to \(\lambda_k \). Then, for any function \(f \in L^2(M) \), 
\[
\pi_k(f) = \sum_{j=1}^{d_k} \langle f, \phi_{k_j} \rangle _{L^2(M)} \phi_{k_j},
\]  
where, \( \{\phi_{k_j} \}_{j=1}^{d_k} \) forms an orthonormal basis for \( \ker(-\Delta_g - \lambda_k) \).  

\medskip
Now for $m> 1$ being fixed, we consider the operator $\mathcal{L}_g$ (cf. \eqref{log-Laplace}) on $(M,g)$, 
which is an unbounded operator on \(L^2(M) \), with the domain of definition 
$$\mathcal{D}(\mathcal{L}_g) :=\left\{ u \in L^2(M) \,\middle|\,\, \sum_{k=0}^\infty |\lambda_k +m|^2|\ln(\lambda_k+m)|^2 |\langle u, \phi_k \rangle|^2 < \infty \right\}.$$ 
For any $s>0$, the inequality 
\begin{equation}\label{log_s}
|\ln(\lambda_k+m)|^2 <(\lambda_k+m)^{2s},\quad k=\{0\}\cup\mathbb{N},\quad\mbox{and }m>1,\end{equation}
yields the Sobolev space\footnote{The Sobolev space $H^\alpha(M)$,
 $\alpha\in\mathbb{R}$ can be considered as defined through equivalence of norms as $c_\alpha\|(\mathbb{I}-\Delta_g)^\frac{\alpha}{2} v\|_{L^2(M)}\leq \|v\|_{H^\alpha(M)} \leq \widetilde{c}_\alpha \|(\mathbb{I}-\Delta_g)^\frac{\alpha}{2} v\|_{L^2(M)}$ for $c_\alpha, \widetilde{c}_\alpha>0$, and $v\in H^\alpha(M)$.}

\begin{equation}\label{H(M)}
H^{2+2s}(M)\subset \mathcal{D}(\mathcal{L}_g):=\mathbb{H}(M).
\end{equation}
We now define the  operator $\mathcal{L}_g$ through its action over its domain of definition $\mathbb{H}(M)$, as
\begin{equation}\begin{aligned}\label{L_g}
    \mathcal{L}_g u &=(-\Delta_g+m\mathbb{I})\circ   \ln [(-\Delta_g+m\mathbb{I})] u ,\quad u\in \mathbb{H}(M)\\
    &= \sum_{k=0}^\infty  (\lambda_k+m)\ln(\lambda_k+m) \pi_k (u).
\end{aligned}\end{equation}
The symbol and order of the pseudodifferential operator $\mathcal{L}_g$ will be discussed in Section \ref{sym-ord}. 

\medskip
In this regard, let us also define the log-Laplacian operator $\ln [(-\Delta_g+m\mathbb{I})]$. Let us call $\mathcal{A}_g=(-\Delta_g+m\mathbb{I})$.
Here we define the operator $\ln\mathcal{A}_g$  through spectral resolution as
\begin{equation}\label{log_A_g}
   \ln\mathcal{A}_g u =\ln [(-\Delta_g+m\mathbb{I})] u 
    = \sum_{k=0}^\infty  \ln(\lambda_k+m) \pi_k (u),
\end{equation}
with its domain of definition
\begin{equation}\label{dom_log}\mathcal{D}(\ln\mathcal{A}_g) :=\left\{ u \in L^2(M) \,\middle|\,\, \sum_{k=0}^\infty |\ln(\lambda_k+m)|^2 |\langle u, \phi_k \rangle|^2 < \infty \right\},\end{equation}
and thanks to \eqref{log_s}, one finds $H^{2s}(M) \subset \mathcal{D}(\ln\mathcal{A}_g):=\mathcal{H}(M)$ for every $s>0$. 
\subsection*{Direct Problem}
Let $f \in C_0^\infty(\mathcal{O})$ and consider the direct problem 
\begin{equation}\label{equation_u}
 \mathcal{L}_gu=f\quad\mbox{ in }M.   
\end{equation}

\medskip
Let us first observe that $0$ is not an eigenvalue of our operator $\mathcal{L}_g$ defined in \eqref{L_g} for $m>1$.

\medskip
If not, suppose $0$ is an eigenvalue of the operator $\mathcal{L}_g$. Then there exists a nonzero function $u\in \mathbb{H}(M)$ such that
\begin{align*}
   \langle u,\, \mathcal{L}_g\, u \rangle = 0.
\end{align*}
Expanding $u$ in terms of the orthogonal $L^2(M)$ eigenbasis $\{\phi_k\}_{k=1}^\infty$, i.e. $u=\sum_{k=1}^\infty \langle \phi_k, u \rangle \phi_k$, we obtain
\begin{align*}
   \left\langle \sum_{k=0}^\infty \langle \phi_k, u \rangle \phi_k,\ \sum_{k=0}^\infty (\lambda_k + m)\ln(\lambda_k + m)\, \langle \phi_k, u \rangle \phi_k \right\rangle = 0,
\end{align*}
which implies
\begin{align*}
   \sum_{k=0}^\infty |\langle \phi_k, u \rangle|^2\, (\lambda_k + m)\ln(\lambda_k + m) = 0.
\end{align*}
Since each term in the sum is non-negative, thanks to the fact $\lambda_k\geq 0$ for all $k \geq 0$,  and $m>1$, it follows that
\begin{align*}
   \langle \phi_k, u \rangle = 0 \quad \text{for all } k \geq 0,
\end{align*}
which implies $u = 0$. This contradicts our assumption that $u$ is nonzero.

\medskip
Therefore, the (linear) operator $\mathcal{L}_g$ is injective on its domain of definition. Further, let us define its inverse acting on $L^2(M)$ spectrally as:
\[\left( \mathcal{L}_g \right)^{-1} f = \sum_{k=0}^\infty \frac{1}{(\lambda_k + m)\ln(\lambda_k + m)} \langle \phi_k, f \rangle \phi_k, \quad f\in L^2(M),
\]
and one notes that, 
$$\mathcal{L}_g \circ \left( \mathcal{L}_g \right)^{-1}f= \sum_{k=0}^\infty \langle \phi_k, f \rangle \phi_k=f, \quad f\in L^2(M),$$
hence by considering 
$$ u= \left( \mathcal{L}_g \right)^{-1}f$$
one finds the unique solution of \eqref{equation_u} in $\mathbb{H}(M)\subset L^2(M)$, for $f\in C^\infty_0(\mathcal{O})$ there.

 \subsection*{Cauchy Data}
Let $u\in \mathbb{H}(M)$ solves \eqref{equation_u} for some $f\in C^\infty_0(\mathcal{O})$, and we define the corresponding Cauchy data set as
\[ \mathcal{C}^{\mathcal{O}}_{M, g}= \{(u|_\mathcal{O},  \mathcal{L}_g u|_\mathcal{O}), \,\,\, \mathcal{L}_g u = 0 \,\mbox{ in } \, M \setminus \overline{\mathcal{O}} \}.\]  This dataset contains solutions to the equation \(\mathcal{L}_gu = 0 \) and serves as the basis on the anisotropic Calder\'{o}n problem for the operator $\mathcal{L}_g$. 

\subsection*{Inverse Problem}
We want to study the inverse problem whether the Cauchy data $\mathcal{C}^{\mathcal{O}}_{M, g}$ determines $(M, g)$ or not, and we have the following to claim as our main result of this paper.
 \begin{theorem}\label{th1}
Let \( (M_i, g_i), \, i = 1,2 \) be two closed, connected Riemannian manifolds with \( \dim(M_i) > 2 \). Assume that \(\mathcal{O} \subset M_1 \cap M_2 \) is a non-empty open subset where the metrics agree, i.e., $(\mathcal{O}, g_1) = (\mathcal{O}, g_2):= (\mathcal{O}, g)$. The equality of the Cauchy data set over $(\mathcal{O}, g)$, i.e. $\mathcal{C}^\mathcal{O}_{M_1, g_1} = \mathcal{C}^\mathcal{O}_{M_2, g_2}$, implies the existence of a diffeomorphism $\Phi: M_1 \to M_2$ such that $\Phi^\ast g_2 = g_1$.
\end{theorem}

\medskip
 
\begin{remark}
In recent years, the research of non-local inverse problems has been quite dynamic. Starting from the article \cite{GSU20}, without being exhaustive we mention some
references \cite{GLX17, RS20, RS2020, GRSU20, MLR20, Cov20, Li20, HL20, BGU21, KLW22, CGR22, Gho22, Zim23, HU24, HLW24, Das25}. However, the operator's order (cf. Section \ref{sym-ord}) in consideration falls within $(0, 2)$. This study extends that limit to $2+$ using this example of $(-\Delta_g+m\mathbb{I})\circ   \ln [(-\Delta_g+m\mathbb{I})]$, where $m>1$.
\end{remark}
In general, these nonlocal problems are motivated by many models, such as diffusion process \cite{AMRT10}, image processing \cite{GO08}, stochastic process \cite{BV16} and so on.

\medskip
The paper is organized as follows. In Section \ref{sec2}, we introduce the fractional power of the elliptic operator using the semigroup approach. Further, we emphasize the symbol and order of these non-local operators.  Section \ref{sec3} is devoted to the proof of our main result.

\section{Functional Calculus}\label{sec2}
\subsection{Heat semigroup approach in defining the  operator log-Laplacian}
Let us recall $\mathcal{A}_g=(-\Delta_g+m\mathbb{I})$.
Here, we define the operator $\ln\,\mathcal{A}_g$ using the semigroup technique, which is equivalent to the prior spectral definition \eqref{log_A_g}.

\medskip
Let \(e^{-t\mathcal{A}_g} \) denote the associated heat semigroup on \(L^2(M)\). Let $v\in L^2(M)$, we define \begin{equation}\label{hk1} e^{-t\mathcal{A}_g} v(x) = \int_M \mathcal{K}_{\mathcal{A}_g}(t,x, y) v(y) \, dV_g(y), 
\end{equation}
where $\mathcal{K}_{\mathcal{A}_g}(t,x,y)$ is the heat kernel for the heat semigroup $e^{-t\mathcal{A}_g}$, which is defined as 
\begin{align}
    \mathcal{K}_{\mathcal{A}_g}(t,x,y) &= \sum_{k=0}^\infty e^{-t(\lambda_k+m)} \phi_k(x)\phi_k(y) \\
     &= e^{-mt}\sum_{k=0}^\infty e^{-t\lambda_k} \phi_k(x)\phi_k(y)=e^{-mt}\mathcal{K}_{\Delta_g}(t,x,y),\label{2.1}
\end{align}
here $\mathcal{K}_{\Delta_g}$ is the heat kernel associated with the heat semigroup $e^{t\Delta_g}$. Grigor’yan, A. \cite{grigor97} provides the following pointwise estimate of the heat kernel $\mathcal{K}_{\Delta_g}(t,x,y)$.
\begin{theorem}[\cite{grigor97}]\label{th2.1}
Let $x, y$ be two points on an arbitrary smooth connected compact Riemannian manifold $M$, and let $t \in (0,\infty)$. Then 
\begin{equation}
   \left| \mathcal{K}_{\Delta_g}(t, x, y) \right| \leq \frac{C}{t^{n/2}} \,\,e^{ -\frac{c\,d^2_g(x, y)}{t} }
\end{equation}

where $C>0$, and \( d_g(x, y) \) is the Riemannian distance between \( x \) and \( y \).
\end{theorem}

\noindent
 
\medskip
Let $x, y$ be two points on an arbitrary smooth connected compact Riemannian manifold $M$, and let $t \in (0,\infty)$. Then \eqref{2.1}, and Theorem\,\ref{th2.1}  yield
\begin{equation}
   \left| \mathcal{K}_{\mathcal{A}_g}(t, x, y) \right| \leq \frac{C}{t^{n/2}} \,e^{-tm}\, e^{ -\frac{c\,d^2_g(x, y)}{t}},\label{ges1}
\end{equation}
where $C>0$, and \( d_g(x, y) \) is the Riemannian distance between \( x \) and \( y \).

\medskip
Since,\footnote{Let $I(\lambda)=\int_0^\infty \frac{(e^{-t}-e^{-t\lambda})}{t}\, dt$, for $\lambda>0$. By doing differentiation under the integral sign \cite[pp. 237]{Rudin64}, we find  $I^\prime(\lambda)=\frac{1}{\lambda}$. Thus $I(\lambda)=\ln\lambda +c$. Thanks to $I(1)=0$, we find $c=0$ here, hence we conclude $I(\lambda)=\ln\lambda$.}
\begin{equation}
    \ln\,\lambda=\int_0^\infty \frac{(e^{-t}-e^{-t\lambda})}{t}\, dt,\quad \lambda>0.
\end{equation}

\medskip
Using the functional calculus (\cite[Ch. 31]{LAX02}), we now define the operator
\begin{equation}
    \ln\,\mathcal{A}_g:=\int_0^\infty \frac{(e^{-t}\mathbb{I}-e^{-t\mathcal{A}_g})}{t}\, dt,
\end{equation}
and its action over its domain of definition as 
\begin{equation}\label{log_A_g_s}
    \ln\,\mathcal{A}_gv(x):=\int_0^\infty \frac{(e^{-t}\mathbb{I}-e^{-t\mathcal{A}_g})v(x)}{t}\, dt,
\end{equation}

where the integral converges in $L^2(M)$ for $v\in \mathcal{H}(M) ,\,\text{where}\,\,\mathcal{D}(\ln\mathcal{A}_g)=\mathcal{H}(M)$ (cf. \eqref{dom_log}), which we are proving next.

\subsection*{Well-definedness of \eqref{log_A_g_s}}
Let us begin with recalling the Gaussian-type upper bound \eqref{ges1} for the heat kernel, and by setting \texorpdfstring{$H_M(z) = e^{-c_1 z^2}$}{HM(z) = exp(-c1 z^2)}\,\,\footnote{Here $z$ stands for the mapping $z: M\times M \times (0,\infty)\to [0,\infty)$ which given by $(x,y,t) \mapsto \frac{d_g(x,y)}{\sqrt{t}}$, where \( d_g \) is the Riemannian  metric on \(( M,g )\).}, we have 
\begin{align}
|e^{-t\mathcal{A}_g} v(x)| &\leq C e^{-tm}\,\,t^{-n/2} \int_M H_M\left( \frac{d_g(x,y)}{\sqrt{t}} \right) |v(y)| \, dV(y) \notag\\
&\,\,= Ce^{-tm}\, \int_M H_M(z) |v(y)| \, dV(z) \notag\\
&\leq C e^{-t}\,\|H_M\|_{L^2(M)} \|v\|_{L^2(M)}, \quad(\mbox{as }m>1),\label{4w}
\end{align} 
Let us now write
\begin{equation}\label{wdl}\begin{aligned}
    \int_0^\infty \frac{(e^{-t}\mathbb{I}-e^{-t\mathcal{A}_g})v(x)}{t}\, dt 
    &=\int_0^1 \frac{(e^{-t}\mathbb{I}-e^{-t\mathcal{A}_g})v(x)}{t}\, dt \\
    &\quad +\int_1^\infty \frac{(e^{-t}\mathbb{I}-e^{-t\mathcal{A}_g})v(x)}{t}\, dt.
\end{aligned}\end{equation}
Let us justify below that the first integral 
\begin{align}
\int_0^1 \frac{(e^{-t}\mathbb{I}-e^{-t\mathcal{A}_g})v(x)}{t}\, dt &= \int_0^1 \partial_t (v_0-v)(\theta(t),x))\,dt
\end{align}
where, $v(t,x)=e^{-t\mathcal{A}_g}v(x)$, and $v_0(t,x)=e^{-t}v(x)$ in $(0, \infty)\times M$; and $\theta(t)\in (0,t)$ appears as some intermediate point due to applying the classical mean-value theorem.

Thus 
\begin{align}
\int_0^1 \frac{(e^{-t}\mathbb{I}-e^{-t\mathcal{A}_g})v(x)}{t}\, dt &=\int_0^1 \left(-e^{-t}v(x)+e^{-t\mathcal{A}_g}\mathcal{A}_gv(x)\right)(\theta(t), x)\,dt\notag\\
&\leq \left(|v(x)| + C\|H_M\|_{L^2(M)}\|\mathcal{A}_gv\|_{L^2(M)}\right)\int_0^1 e^{-\theta(t)}\,dt\notag\\
&\leq \left(|v(x)|+ C\|H_M\|_{L^2(M)}\|\mathcal{A}_gv\|_{L^2(M)}\right),
\end{align}
where on the second-to-last line we have used the estimate \eqref{4w}. 

\medskip
Next, we consider the second integral in the r.h.s. of \eqref{wdl}; and using \eqref{4w}, we simply find 
\begin{equation*}
     \left|\int_1^\infty \frac{(e^{-t}\mathbb{I}-e^{-t\mathcal{A}_g})v(x)}{t}\, dt \right| \leq (1+C \|H_M\|_{L^2(M)})\,\| v\|_{L^2(M)}\int_1^\infty\, \frac{e^{-t}}{t}\, dt <\infty.
\end{equation*}
Hence, the integral \eqref{log_A_g_s} converges pointwise for smooth $v\in C^\infty(M)$, to have 
\begin{equation}\label{esti}
|\ln\,\mathcal{A}_gv(x)|\leq \left(|v(x)|+ C\left(\|v\|_{L^2(M)}+\|\mathcal{A}_gv\|_{L^2(M)}\right)\right), \quad x\in M.
\end{equation}
Since $(M,g)$ is compact, this also gives the $L^2(M)$-ness of $\ln\,\mathcal{A}_gv$ for $v\in H^2(M)$, as 
\begin{equation}\label{esti2}
\|\ln\,\mathcal{A}_gv\|_{L^2(M)}\leq C\left(\|v\|_{L^2(M)}+ \|\mathcal{A}_gv\|_{L^2(M)}\right).
\end{equation}
\medskip
Following this (cf. \eqref{esti}-\eqref{esti2}), we define 
\begin{align}
(-\Delta_g+m\mathbb{I}) \circ   \ln [(-\Delta_g+m\mathbb{I})]v(x)&:= \mathcal{A}_g\circ\ln (\mathcal{A}_g) \,v(x)\notag\\
&=\ln (\mathcal{A}_g)\circ \mathcal{A}_g \,v(x)\notag\\
&=   \int_0^\infty \frac{\left(e^{-t}\mathbb{I}-e^{-t\mathcal{A}_g}\right) \mathcal{A}_gv(x)}{t} \, dt,\label{df_h}
\end{align}
where the above integral converges in $L^2(M)$, for $v\in H^4(M)$.

\subsection{Symbol and order}\label{sym-ord}
We recall the following integral representation. \footnote{It essentially follows from this identity $\lambda^{-a}=\frac{1}{\Gamma(a)}\,\int_0^\infty \frac{e^{-t\lambda}}{t^{1-a}}\,dt$, for $a\in (0,1)$ and $\lambda>0$. }
  
\begin{equation}\label{lambda}
\lambda^s =\frac{1}{\Gamma(-s)}\,\int_0^\infty \frac{(e^{-t\lambda}-1)}{t^{s+1}}\,dt,\quad s\in (0,1), \quad\lambda\geq 0.
\end{equation}
where $\Gamma(-s)=\frac{\Gamma(1-s)}{-s}$. 

\medskip
Let $v\in C_0^\infty(M)$. Then the functional calculus or the semigroup representation of the fractional power of $\mathcal{A}_g=(-\Delta_g+m\mathbb{I})$ becomes
\begin{equation}\label{frac_def}
(-\Delta_g+m\mathbb{I})^sv(x)=\mathcal{A}_g^sv(x) := \frac{1}{\Gamma(-s)} \int_0^\infty \frac{\left(e^{-t\mathcal{A}_g}-\mathbb{I}\right)v(x)}{t^{1+s}} \, dt,\quad s\in (0,1).
\end{equation}
 We begin with the following result
\begin{theorem}[\cite{FGKRSU25}]
  Let  $\beta > 0$ and $g$ be a $C^\infty$ Riemannian metric on $\mathbb{R}^n$, $n \geq 2$, which agrees with the Euclidean metric outside of a compact set. Then the operator $(-\Delta_g)^\beta\in \Psi^{2\beta}_{1,0}(\mathbb{R}^n)$, and is a classical elliptic pseudodifferential operator on $\mathbb{R}^n$, with principal symbol 
\[
\left(\sum_{j,k=1}^n g^{jk}(x)\xi_j\xi_k\right)^\beta, \quad (x,\xi) \in T^*\mathbb{R}^n \setminus \{0\}.
\] 
\end{theorem}
\begin{remark}
In the preceding finding, we understand $(-\Delta_g)^\beta$ as follows. Let $\beta=m+\alpha$, where $m\in\mathbb{Z}$ and $\alpha\in (0,1)$, then $(-\Delta_g)^\beta=(-\Delta_g)^m\circ (-\Delta_g)^\alpha$.

\end{remark}
Next, we review the symbol class definitions.

\medskip
Assume \(m, \rho, \delta \in \mathbb{R} \) and \(0 \leq \rho, \delta \leq 1 \). We define the symbol class \( S^m_{\rho,\delta}(\mathbb{R}^n \times \mathbb{R}^n \times \mathbb{R}^n) \) to consist of all functions \( p \in C^\infty(\mathbb{R}^n \times \mathbb{R}^n \times \mathbb{R}^n) \) such that, for every compact set \( K \subset \mathbb{R}^n \) and for all multi-indices \( \alpha, \beta \), there exists a constant \( C_{K,\alpha,\beta} > 0 \) such that
\[
|\partial_x^\alpha \partial_\xi^\beta p(x, \xi)| \leq C_{K,\alpha,\beta} (1 + |\xi|)^{m - \rho |\beta| + \delta |\alpha|}, \quad \text{for all } x \in K,\, \xi \in \mathbb{R}^n.
\]

\medskip
In particular, the subclass \( S^0_{1,0}(\mathbb{R}^n) \) consists of symbols \( p(x, \xi) \in C^\infty(\mathbb{R}^n \times \mathbb{R}^n \times \mathbb{R}^n)) \) for which the estimates
\[
|\partial_x^\alpha \partial_\xi^\beta p(x, \xi)| \leq C_{K,\alpha,\beta} (1 + |\xi|)^{ - |\beta|}
\]
hold uniformly for all \( x \in K  \), where \( K \) is any compact set.

 The space \( \Psi^0_{1,0}(\mathbb{R}^n) \) of pseudodifferential operators on \( \mathbb{R}^n \)
    associated with the symbol class \( S^0_{1,0}(\mathbb{R}^n \times \mathbb{R}^n \times \mathbb{R}^n) \)
consisting of operators of the form
\[
A u(x) = \frac{1}{(2\pi)^n} \int_{\mathbb{R}^n} \int_{\mathbb{R}^n} e^{i(x - y)\cdot \theta} \, a(x, y, \theta) \, u(y) \, dy \, d\theta, \quad u \in C^\infty_0(\mathbb{R}^n),
\]
where the amplitude \( a(x, y, \theta) \in S^0_{1,0}(\mathbb{R}^n \times \mathbb{R}^n \times \mathbb{R}^n) \). For more details, we refer to  \cite {GS94}.

\subsection*{Symbol class of the logarithmic Laplacian:}
The asymptotic expansion of $\lambda^s$ (cf. \eqref{lambda}) for $s\to 0+$ stands as
$$ \lambda^s= 1+ s\,\ln\lambda+ \mathcal{O}(s^2),$$
or,
$$\ln\lambda =\left.\frac{d}{ds}\right|_{s=0}\lambda^s,\quad \lambda>0.$$
This infers that, thanks to functional calculus, the logarithmic Laplacian $\ln \mathcal{A}_g$ which has been defined in \eqref{log_A_g} or \eqref{log_A_g_s}, enjoy the following representation 
$$\ln \mathcal{A}_g= \left.\frac{d}{ds}\right|_{s=0} (-\Delta_g+m\mathbb{I})^s.$$
This approach has been discussed in \cite{LAWT21, CHW23} to define the near-zero logarithmic Laplacian operator. 

Since the principal symbol of $(-\Delta_g+m\mathbb{I})^s$ is $\left(\sum_{j,k=1}^n g^{jk}(x)\xi_j\xi_k +m\right)^s$, $s\in (0,1)$. Therefore, the principal symbol of $\ln \mathcal{A}_g$ would be
\[
\sigma(\ln \mathcal{A}_g)(x,\xi) = \ln\left(\sum_{j,k=1}^n g^{jk}(x)\xi_j\xi_k +m\right).
\]
Since for any multi-index \( \beta \in \mathbb{N}^n \),  
$$
|\partial_\xi^\beta \ln(|\xi|_g^2 + m)| \leq C_\beta (1 + |\xi|)^{-|\beta|},
$$
where $|\xi|_g^2=\left( \sum_{j,k=1}^n g^{jk}(x)\xi_j\xi_k +m\right)$. 

\medskip
Thus, the logarithmic Laplacian is a classical pseudodifferential operator of order zero:
\[
\ln \mathcal{A}_g \in \Psi^0_{1,0}
\]
with logarithmic principal symbol.

\medskip
Similarly, the principal symbol of the operator $\mathcal{L}_g=(-\Delta_g + m\mathbb{I})\circ\ln\,(-\Delta_g + m\mathbb{I})$ would be
$$\sigma(\mathcal{L}_g)(x, \xi) = (|\xi|_g^2 + m)\,\ln\,(|\xi|_g^2 + m).$$
Since for any \( \varepsilon > 0 \), we have the estimate:
\[
|\xi|^2 \ln |\xi| \leq C_\varepsilon (1 + |\xi|)^{2 + \varepsilon}.
\]
Furthermore, the derivatives of the symbol satisfy:
\[
\partial_x^\alpha \partial_\xi^\beta \sigma(\mathcal{L}_g)(x, \xi) = O\left((1 + |\xi|)^{2 + \varepsilon - |\beta|} \right).
\]
Thus $\sigma(\mathcal{L}_g)(x, \xi)\in S^{2 + \varepsilon}_{1,0}$, i.e, 
$$\mathcal{L}_g=(-\Delta_g + m\mathbb{I})\circ\ln[(-\Delta_g + m\mathbb{I})] \in \Psi^{2 + \varepsilon}_{1,0},  \quad\text{for any } \varepsilon > 0,$$
establishing it is an operator of order $2+$.

\bigskip

\section{Proof of Main Result}\label{sec3}
We recall, $\mathcal{O}\subset M_i$ for $i=1,2$. Let $f\in C^\infty_0(\mathcal{O})$, and we uniquely solve  the equation 
\begin{equation}\label{eq_i}
\mathcal{L}_{g_i} u^f_i=f \quad\mbox{in }(M_i, g_i),
\end{equation}
for $u_i:=u^f_i\in \mathbb{H}(M_i)$. 

\medskip
By our hypothesis, we have  
\begin{equation}\label{cdata} (g_1|_{\mathcal{O}},\,u^f_1|_\mathcal{O},  \mathcal{L}_{g_1} u^f_1|_\mathcal{O})=(g_2|_{\mathcal{O}},\,u^f_2|_\mathcal{O},  \mathcal{L}_{g_2} u^f_2|_\mathcal{O}).
\end{equation}  
\subsection*{Regularity of \texorpdfstring{$u_i$}{u	extunderscore i}} This can be observed that $u^f_i\in C^\infty(M_i)$ for $f\in C^\infty_0(\mathcal{O})$. As it follows 
\begin{equation}\label{Aleq}
\mathcal{L}_{g_i}\left( \mathcal{A}_{g_i} u^f_i\right)=\mathcal{L}_{g_i}\circ \mathcal{A}_{g_i} u^f_i=\mathcal{A}_{g_i}\circ \mathcal{L}_{g_i} u^f_i=\mathcal{A}_{g_i} f \quad\mbox{in }(M_i, g_i)\end{equation} 
where, $\mathcal{A}_{g_i}=(-\Delta_{g_i}+m\mathbb{I})$ in $(M_i, g_i)$ for $i=1,2$. 

By solving \eqref{Aleq}, with the source term   $\mathcal{A}_{g_i}f\in  L^2(M_i)$, we find $\mathcal{A}_{g_i} u^f_i\in  \mathbb{H}(M_i)$. Continuing the process, one shows $\mathcal{A}^k_{g_i} u^f_i\in  \mathbb{H}(M_i)\subset L^2(M_i)$ for every $k\in \mathbb{N}$, hence establishing the claim that $u^f_i\in \underset{l\in\mathbb{N}\setminus\{1,2\}}{\cap}\, H^l(M_i)=\underset{l\in\mathbb{N}\cup\{0\}}{\cap}\, H^l(M_i)$; or in $C^\infty(M_i)$ as the closed infinite intersections of the Sobolev spaces $\{H^l(M_i)\}_{l\in\mathbb{N}\cup \{0\}}$.  

\subsection*{Integral identities}
Further, by realizing 
\[
\mathcal{L}_{g_i}\left( \mathcal{A}^{k-1}_{g_i} u^f_i\right)=\mathcal{L}_{g_i}\circ \mathcal{A}^{k-1}_{g_i} u^f_i=\mathcal{A}^{k-1}_{g_i}\circ \mathcal{L}_{g_i} u^f_i=\mathcal{A}^{k-1}_{g_i} f \quad\mbox{in }(M_i, g_i), \quad k\in \mathbb{N};
\]
and 
\[\mathcal{A}_{g_1}^{k-1}f= \mathcal{A}_{g_2}^{k-1}f\quad\mbox{in }(M_i, g_i), \quad k\in \mathbb{N}; \] 
thanks to our choice $f\in C^\infty_0(\mathcal{O})$, (as $\mathcal{A}_{g_i}$ are the local operators) and there $g_1|_{\mathcal{O}}=g_2|_{\mathcal{O}}$. 

\medskip
Therefore, the equality of the Cauchy data set (cf. \eqref{cdata}) allow us to write: 
\begin{equation}
    \left(\left(\mathcal{A}_{g_1}^{k-1} u_1\right) |_{\mathcal{O}},\,\, \left(\mathcal{L}_{g_1} \mathcal{A}_{g_1}^{k-1} u_1 \right)\big|_{\mathcal{O}}\right) =\left(\left(\mathcal{A}_{g_2}^{k-1} u_2\right) |_{\mathcal{O}},\,\, \left(\mathcal{L}_{g_2} \mathcal{A}_{g_2}^{k-1} u_2 \right)\big|_{\mathcal{O}}\right), \qquad k \in \mathbb{N}.\label{k-cdata}
\end{equation}


Indeed, the regularity of \( u_i \) ensures that \( u_i \in \mathcal{H}(M_i) \). Consequently, employing the semigroup definition of \( \mathcal{L}_{g_i} \) in~\eqref{df_h} together with~\eqref{k-cdata}, we deduce the following integral identities:

\begin{equation}
    \int_0^\infty \left(e^{-t\mathcal{A}_{g_1}}\mathcal{A}_{g_1}^{k} u_1 - e^{-t\mathcal{A}_{g_2}}\mathcal{A}_{g_2}^{k} u_2 \right)(t, x)\, \frac{dt}{t} = 0, \quad\label{int_idn}
\end{equation}
for $x\in \mathcal{O}$, and $ k \in \mathbb{N}$.

\medskip
On the domain
$\mathcal{D}((-\Delta_{g_i}+m\mathbb{I})^j)$, the operators commute: i.e,\,\, $e^{-t\mathcal{A}_{g_i}} \mathcal{A}_{g_i}^j = \mathcal{A}_{g_i}^j e^{-t\mathcal{A}_{g_i}}$ \,\,for all $t \geq 0$. In particular we have
\begin{equation}
    e^{-t\mathcal{A}_{g_i}} (\mathcal{A}_{g_i}^j u_i) =(-1)^j\, \partial_t^j (e^{-t\mathcal{A}_{g_i}} u_i), \qquad j\in\mathbb{N}\,\,\text{and}\,\, (t,x)\in [0,\infty)\times (M_i, g_i).\label{sg}
\end{equation}
Substituting \eqref{sg} into \eqref{int_idn} gives
\begin{equation}
    \int_0^\infty \partial_t^{k} \left(e^{-t\mathcal{A}_{g_1}}u_1 - e^{-t\mathcal{A}_{g_2}} u_2\right)(t, x)\, \frac{dt}{t} = 0,\label{int_idn2}
\end{equation}
for $x \in\mathcal{O}$ and $k\in\mathbb{N}$.

\medskip
We wish to perform the integration by-parts (consecutively $k$ times) on the above integral to write
\begin{equation}\label{int_idn3}
\int_0^\infty \left(e^{-t\mathcal{A}_{g_1}}u_1 - e^{-t\mathcal{A}_{g_2}} u_2 \right)(t, x)\, \frac{dt}{t^{1 + k}} = 0,
\end{equation}
for $x \in \mathcal{O}$, and $k\in\mathbb{N}$.

\medskip
Let us justify, why the above integral \eqref{int_idn3} is well-defined, for $x\in\mathcal{O}$, and for all $k\in\mathbb{N}$. 
\subsection*{Well-definedness of \texorpdfstring{\eqref{int_idn3}}{(3)}}
Let $\omega\Subset\mathcal{O}$ be some arbitrary non-empty open set. We know $g_1|_{\mathcal{O}}=g_2|_{\mathcal{O}}$. The restriction of the function $\left(e^{-t\mathcal{A}_{g_1}}u_1 - e^{-t\mathcal{A}_{g_2}} u_2 \right)(t, x)$  belongs to $C^\infty((0, \infty) \times \overline{\omega})$, with satisfying the heat equation
\begin{align}
    (\partial_t - \mathcal{A}_{g_1}) \left(e^{-t\mathcal{A}_{g_1}}u_1 - e^{-t\mathcal{A}_{g_2}} u_2 \right)= 0 \quad \text{in } (0,\infty)\times\omega.
\end{align}
with the initial data over $\{0\}\times \omega$
\begin{equation}\label{ini1}
\left(e^{-t\mathcal{A}_{g_1}}u_1 - e^{-t\mathcal{A}_{g_2}} u_2\right)(0, x)\big|_{x\in\omega}=(u_1-u_2)(x)\big|_{x\in\omega}=0.
\end{equation}
Moreover, as follows from  \eqref{sg} and \eqref{k-cdata}, that the all right derivatives w.r.t. time variable of $\left(e^{-t\mathcal{A}_{g_1}}u_1 - e^{-t\mathcal{A}_{g_2}} u_2 \right)$, at $t=0$ exists, and equals to $0$, i.e. 
\begin{equation}\label{ini2}
\partial^{j}_t\Big[\left(e^{-t\mathcal{A}_{g_1}}u_1 - e^{-t\mathcal{A}_{g_2}} u_2\right)\Big](0, x)\big|_{x\in\omega}=(\mathcal{A}^j_{g_1}u_1-\mathcal{A}^j_{g_2}u_2)(x)|_{x\in\omega}=0, \qquad j\in\mathbb{N}.
\end{equation}
So, $\left(e^{-t\mathcal{A}_{g_1}}u_1 - e^{-t\mathcal{A}_{g_2}} u_2 \right)(t, x)$  belongs to $C^\infty([0, \infty) \times \overline{\omega})$.

\medskip
Now we will argue as we did in showing the well-definedness of \eqref{log_A_g_s}. We write the integral \eqref{int_idn3} as
\begin{equation}\label{wdl2}\begin{aligned}
    \int_0^\infty \left(e^{-t\mathcal{A}_{g_1}}u_1 - e^{-t\mathcal{A}_{g_2}} u_2 \right)(t, x)\, \frac{dt}{t^{1 + k}}
    &=\int_0^1 \left(e^{-t\mathcal{A}_{g_1}}u_1 - e^{-t\mathcal{A}_{g_2}} u_2 \right)(t, x)\, \frac{dt}{t^{1 + k}} \\
    &\quad +\int_1^\infty \left(e^{-t\mathcal{A}_{g_1}}u_1 - e^{-t\mathcal{A}_{g_2}} u_2 \right)(t, x)\, \frac{dt}{t^{1 + k}}.
\end{aligned}\end{equation}
Let us justify below that the first integral. We will be using the Taylor series expansion for $\phi\in C^{l-1}[0, 1]\cap C^l(0,1)$, $l\in\mathbb{N}$,
    \begin{equation}\label{phi}
    \phi(t)= \sum_{j=0}^{l-1} \frac{t^j}{j!} \partial^j \phi(0) + \frac{t^l}{(m-1)!}\int_0 ^1 (1-\theta)^{l-1}(\partial^l \phi)(\theta t)\, d\theta.\end{equation}
If $(\partial^j )\phi(0) =0 \,\,\text{for}\,\, 0\leq j\leq l-1, \,\, \phi (t)=t^l \psi(t)$ were $\psi \in C^\infty$ and supp$(\psi)=$ supp$(\phi)$.

\medskip
For $x\in\omega$ being fixed, by choosing $\phi(t):=\left(e^{-t\mathcal{A}_{g_1}}u_1 - e^{-t\mathcal{A}_{g_2}} u_2\right)(t, x)$   
\eqref{ini1}-\eqref{ini2}, we indeed have $\left(e^{-t\mathcal{A}_{g_1}}u_1 - e^{-t\mathcal{A}_{g_2}} u_2\right)(t, x)=t^l\psi_l(t, x)$ for some $\psi_l\in L^\infty(0,1)$, $l\in\mathbb{N}$. This shows the first integral is well-defined:
 \begin{equation}\label{1wed}
     \left|\int_0^1 \left(e^{-t\mathcal{A}_{g_1}}u_1 - e^{-t\mathcal{A}_{g_2}} u_2 \right)(t, x)\, \frac{dt}{t^{1 + k}}\right| <\infty, \qquad k\in\mathbb{N}.
 \end{equation}
Next, we will be considering the second integral in the r.h.s. of  \eqref{wdl2}. Let us recall that we showed in \eqref{4w}:
\begin{equation}\label{recall}
|e^{-t\mathcal{A}_g} v(x)| \leq  C e^{-t}\,\|H_M\|_{L^2(M)} \|v\|_{L^2(M)}, 
\end{equation}
where \(H_M(z) = e^{-c_1 z^2}\). Thus 
\begin{align*}
     &\left|\int_1^\infty \left(e^{-t\mathcal{A}_{g_1}}u_1 - e^{-t\mathcal{A}_{g_2}} u_2 \right)(t, x)\, \frac{dt}{t^{1 + k}} \right|\\
     &\qquad\leq C  \big(\|u_1\|_{L^2(M_1)}\,\|H_{M_1}\|_{L^2(M_1)} +  \| u_2\|_{L^2(M_2)}\, \|H_{M_2}\|_{L^2(M_2)}\big)\int_1^\infty\, \frac{e^{-t}}{t}\, dt<\infty.
\end{align*}
This completes the well-definedness of \eqref{int_idn3}.

\medskip
We now justify the absence of endpoint contributions in \eqref{int_idn3} when integration by parts is performed to \eqref{int_idn2}. At both ends, the contributions are $0$, as it is evident from \eqref{phi} for $t=0$, and \eqref{recall} for $t=\infty$ cases respectively. 

\subsection*{Analysis based on the integral identity \texorpdfstring{\eqref{int_idn3}}{(3)}}
Next, we would like to conclude from the integral identity \eqref{int_idn3}, \begin{equation}\label{claim}
    \left(e^{-t\mathcal{A}_{g_1}}u_1 - e^{-t\mathcal{A}_{g_2}} u_2 \right)(t, x)=0\qquad\mbox{ over }[0, \infty) \times \overline{\omega}.
\end{equation}
Let us fix $x\in\overline{\omega}\subset\mathcal{O}$, and call $\phi(t)=    \left(e^{-t\mathcal{A}_{g_1}}u_1 - e^{-t\mathcal{A}_{g_2}} u_2 \right)$. As evident from  \eqref{recall}, 
\[t\in (0,\infty)\longmapsto \phi(t)\in L^2(0,\infty).\]
Let us also see from \eqref{sg} that,
\[\phi^\prime(t)=\frac{d}{dt}\phi(t)= \left(e^{-t\mathcal{A}_{g_1}}\mathcal{A}_{g_1}u_1 - e^{-t\mathcal{A}_{g_2}}\mathcal{A}_{g_2}u_2 \right),\]
and using \eqref{recall}, we have
\[t\in (0,\infty)\longmapsto \phi^\prime(t)\in L^2(0,\infty),\]
since
\begin{align*}\int_0^\infty|\phi^\prime(t)|^2\,dt &=\int_0^\infty \big|e^{-t\mathcal{A}_{g_1}}\mathcal{A}_{g_1} u_1- e^{-t\mathcal{A}_{g_2}}\mathcal{A}_{g_2}u_2\big|^2\, dt\\ 
&\leq  C  \big(\|\mathcal{A}_{g_1}u_1\|_{L^2(M_1)}\, \|H_{M_1}\|_{L^2(M_1)}+\|\mathcal{A}_{g_2}u_2\|_{L^2(M_2)}\,\|H_{M_2}\|_{L^2(M_2)}\big)\,  \int_0^\infty e^{-2t}<\infty.
\end{align*}
\subsection*{Hardy's inequality and application}
Let us apply the classical Hardy's inequality on $\phi$, \footnote{The celebrated Hardy inequality states that, if $1 < p < \infty$ and if $u$ is a locally
absolutely continuous function on $(0, \infty)$ with $\liminf_{r\to 0}|u(r)| = 0$, then 
\[\int_0^\infty \frac{|u(r)|^p}{r^p}\, dr\leq (\frac{p}{p-1})^p\,\int_0^\infty |u^\prime(r)|^p\, dr.\]
The constant on the right side is the best possible. 
We refer to \cite{FLW} and reference therein.} to have
$
\int_0^\infty \frac{|\phi(t)|^2}{t^2}\, dt\leq \int_0^\infty |\phi^\prime(t)|^2\, dt.
$
As an application of this, we now claim 
\begin{equation}\label{phi-psi}
    t\in (0,\infty)\longmapsto \psi(t):=\phi(\frac{1}{t})\in L^2(0,\infty),
\end{equation}
as 
\[
\int_0^\infty |\psi(t)|^2\, dt=\int_0^\infty \frac{|\phi(t)|^2}{t^2}\, dt<\infty.
\]
\subsection*{Paley-Wiener theorem and application}
Hence, by considering the Fourier transform of $\psi\in L^2(0,\infty)$ over the upper-half complex plane $H_{+}=\{\xi+\mathrm{i}\eta:\,\,\xi\in\mathbb{R}, \mbox{ and }\eta>0\}$, we conclude
\[
(\xi+\mathrm{i}\eta)\in H_{+}\longmapsto    \widehat{\psi}(\xi+\mathrm{i}\eta)=\int_0^\infty e^{\mathrm{i}\tau(\xi+\mathrm{i}\eta)}\psi(\tau)\,d\tau 
\]
is a holomorphic function, thanks to the classical Paley-Wiener theorem\footnote{
By the classical Paley–Wiener theorem \cite{PW87}: Given any $g\in L^2(0,\infty)$, defining $f$ as
$$H_{+}\ni (x+\mathrm{i}y) \longmapsto f(x+\mathrm{i}y):=\frac{1}{2\pi}\int_0^\infty e^{\mathrm{i}\tau(x+\mathrm{i}y)}g(\tau)\,d\tau,$$ becomes a holomorphic function over upper-half complex plane $H_{+}:=\{x+\mathrm{i}y:\,x\in\mathbb{R},\,\mbox{and}\,y>0\}$, with satisfying   
$$
\sup_{y>0} \int_{-\infty}^{\infty} |f(x+iy)|^2\,dx= \frac{1}{2\pi}\|g\|^2_{L^2(0,\infty)}.
$$
}.

Now by writing the Taylor series expansion of the holomorphic function $\widehat{\psi}(\xi+\mathrm{i}\eta)$, within its radius of convergence at any given point $(\xi+\mathrm{i}\eta)\in H_{+}$, we find
\begin{align*}
\widehat{\psi}(\xi+\mathrm{i}\eta)&=\sum_{j=0}^\infty\left(\int_0^\infty \tau^j\psi(\tau)\,d\tau\right)\frac{(\xi+\mathrm{i}\eta)^j}{j!}\\\notag
&=\sum_{k=1}^\infty\left(\int_0^\infty \frac{\phi(t)}{t^{1+k}}\,dt\right)\frac{(\xi+\mathrm{i}\eta)^{k-1}}{(k-1)!}\\\notag
&=0,
\end{align*}
due to the integral identity \eqref{int_idn3} for all $k\in\mathbb{N}$.

Therefore, we conclude $\psi\equiv0$, consequently $\phi\equiv 0$. This completes the proof of the claim \eqref{claim}.

\medskip
Since $\omega\Subset \mathcal{O}$ was arbitrary, from \eqref{claim}, we finally conclude 
\begin{equation}\label{eql_O}
     (e^{-t\mathcal{A}_{g_1}}\,u_1 - e^{-t\mathcal{A}_{g_2}}\,u_2)\,(t,x)   = 0 \quad \forall \,t > 0, \, x \in \mathcal{O}.
\end{equation}
\subsection*{Equality of the heat kernels}
Let $\mathcal{K}_{\mathcal{A}_{g_i}}(t, x, y)$ denotes the corresponding heat kernels for $i=1, 2$. We will show
\begin{equation}\label{hk}
    \mathcal{K}_{\mathcal{A}_{g_1}}(t, x, y) = \mathcal{K}_{\mathcal{A}_{g_2}}(t, x, y) \quad \forall x, y \in \mathcal{O}, \,\, t > 0.
\end{equation}
Since, over the domain of definition, thanks to functional calculus, we could write 
\[ e^{-t\mathcal{A}_{g_i}} \circ \mathcal{L}_{g_i} = \mathcal{L}_{g_1}\circ e^{-t\mathcal{A}_{g_i}},\quad i=1, 2.\]
Now by fixing some $t=t_0 > 0$, we will be considering the following difference, over $x\in \mathcal{O}$, to find
\begin{align}
 e^{-t_0\mathcal{A}_{g_1}} \circ \mathcal{L}_{g_1}u_1- e^{-t_0\mathcal{A}_{g_2}} \circ \mathcal{L}_{g_2}u_2 &= \mathcal{L}_{g_1}(e^{-t_0\mathcal{A}_{g_1}}u_1)-\mathcal{L}_{g_2}(e^{-t_0\mathcal{A}_{g_2}})\notag\\[2pt]
 &=\int_0^\infty \frac{\left(e^{-t}\mathbb{I}- e^{-t\mathcal{A}_{g_1}}\right)(e^{-t_0\mathcal{A}_{g_1}} u_1)}{t}\,dt \notag\\
 &\qquad-\int_0^\infty \frac{\left(e^{-t}\mathbb{I}- e^{-t\mathcal{A}_{g_2}}\right)(e^{-t_0\mathcal{A}_{g_2}} u_2)}{t}\,dt\quad\mbox{(cf. \eqref{log_A_g_s})}\notag\\[2pt]
&= \int_0^\infty \frac{ e^{-t}\left(e^{-t_0\mathcal{A}_{g_1}}\, u_1 - e^{-t_0\mathcal{A}_{g_2}} \,u_2\right)}{t}\,dt\notag\\
&\qquad-\int_0^\infty \frac{\left( e^{-(t+t_0)\mathcal{A}_{g_1}} u_1 - e^{-(t+t_0)\mathcal{A}_{g_2}} u_2 \right) }{t}\,dt\notag\\[2pt]
&=0, \quad\mbox{(thanks to \eqref{eql_O}).}\label{idn_1}
\end{align}
Since, $\mathcal{L}_{g_i} u_i =f$ in $(M_i, g_i)$ for $i=1, 2$ for $f\in C^\infty_c(\mathcal{O})$ (see \eqref{eq_i}).  
Therefore, from above \eqref{idn_1}, we conclude that 
\begin{equation}
0= e^{-t_0\mathcal{A}_{g_1}} (\mathcal{L}_{g_1} u_1) - e^{-t_0\mathcal{A}_{g_2}} (\mathcal{L}_{g_2} u_2) = e^{-t_0\mathcal{A}_{g_1}}f - e^{-t_0\mathcal{A}_{g_2}}f
\end{equation}
over $x\in \mathcal{O}$.

\medskip
Since $t_0>0$ is arbitrary, thus we have 
\begin{equation}
[e^{-t\mathcal{A}_{g_1}} f](t,x) = [e^{-t\mathcal{A}_{g_2}} f](t,x), \quad \forall t>0,\,\, x\in \mathcal{O}, \qquad\forall f\in C_0^\infty(\mathcal{O}).
\end{equation}
This further implies the equality of the heat kernels, i.e.
\begin{equation}
\mathcal{K}_{\mathcal{A}_{g_1}}(t, x, y) = \mathcal{K}_{\mathcal{A}_{g_2}}(t, x, y) \quad \forall x, y \in \mathcal{O}, \, t > 0,
\end{equation}
and proves our claim \eqref{hk}. 

Since $\mathcal{K}_{\mathcal{A}_{g_i}}(t, x, y)=e^{-mt}\mathcal{K}_{\Delta_{g_i}}(t, x, y)$ for $i=1,2$, $m>0$. Consequently, we the equality of heat kernels of Laplace-Beltrami operators follows as
\begin{equation}
\mathcal{K}_{\Delta_{g_1}}(t, x, y)=\mathcal{K}_{\Delta_{g_2}}(t, x, y)\quad \forall x, y \in \mathcal{O}, \, t > 0.
\end{equation}

\noindent
Finally, we use the following result from \cite[Theorem 1.5]{FGKU25} to finish the proof of our Theorem \ref{th1}.
\begin{theorem}[\cite{FGKU25}]
Let $(N_1, g_1)$ and $(N_2, g_2)$ be smooth connected
complete Riemannian manifolds of dimension $n\geq 2$ without boundary. Let $\mathcal{O}_j \subset N_j$, $j = 1, 2,$ be open nonempty sets. Assume that
$\mathcal{O}_1 =\mathcal{O}_2 := \mathcal{O}$. Assume furthermore that
$$\mathcal{K}_{\Delta_{g_1}}(t, x, y) = \mathcal{K}_{\Delta_{g_2}}(t, x, y),\quad\forall t>0,\,\,\mbox{and }x,y\in\mathcal{O}.$$
Then there exists a diffeomorphism $\varphi : N_1 \mapsto N_2$ such that $\varphi^\ast g_2 =g_1$ on $N_1$.
\end{theorem} 
\medskip
This concludes our discussion of the proof of the Theorem \ref{th1}.
\qed
\qed

\subsection*{Acknowledgment:}
  I am grateful to Prof. Tuhin Ghosh (Harish-Chandra Research Institute, Prayagraj, India)  for his guidance and advice,  all of which helped me to accomplish this project. I'd also want to thank the Department of Atomic Energy, Government of India, for supporting this research through its Fellowship.


\begin{thebibliography}{10}

\bibitem{AMRT10}
{\sc Andreu-Vaillo, F., Maz\'{o}n, J.~M., Rossi, J.~D., and Toledo-Melero, J.~J.}
\newblock {\em Nonlocal diffusion problems}, vol.~165 of {\em Mathematical Surveys and Monographs}.
\newblock American Mathematical Society, Providence, RI; Real Sociedad Matem\'{a}tica Espa\~{n}ola, Madrid, 2010.

\bibitem{KL06}
{\sc Astala, K., and P\"aiv\"arinta, L.}
\newblock Calder\'on's inverse conductivity problem in the plane.
\newblock {\em Ann. of Math. (2) 163}, 1 (2006), 265--299.

\bibitem{BGU21}
{\sc Bhattacharyya, S., Ghosh, T., and Uhlmann, G.}
\newblock Inverse problems for the fractional-{L}aplacian with lower order non-local perturbations.
\newblock {\em Trans. Amer. Math. Soc. 374}, 5 (2021), 3053--3075.

\bibitem{BV16}
{\sc Bucur, C., and Valdinoci, E.}
\newblock {\em Nonlocal diffusion and applications}, vol.~20 of {\em Lecture Notes of the Unione Matematica Italiana}.
\newblock Springer, [Cham]; Unione Matematica Italiana, Bologna, 2016.

\bibitem{Buk08}
{\sc Bukhgeim, A.~L.}
\newblock Recovering a potential from {C}auchy data in the two-dimensional case.
\newblock {\em J. Inverse Ill-Posed Probl. 16}, 1 (2008), 19--33.

\bibitem{MLR20}
{\sc Ceki\'{c}, M., Lin, Y.-H., and R\"{u}land, A.}
\newblock The {C}alder\'{o}n problem for the fractional {S}chr\"{o}dinger equation with drift.
\newblock {\em Calc. Var. Partial Differential Equations 59}, 3 (2020), Paper No. 91, 46.

\bibitem{CHW23}
{\sc Chen, H., Hauer, D., and Weth, T.}
\newblock {A}n extension problem for the logarithmic {L}aplacian.
\newblock {\em https://arxiv.org/pdf/2312.15689\/} (2025).

\bibitem{CR25}
{\sc Chen, X., and R\"{u}land, A.}
\newblock {B}oundary reconstruction for the anisotropic fractional {C}alder{\'{o}}n problem.
\newblock {\em arXiv:2502.12287\/} (2025).

\bibitem{CO24}
{\sc Choulli, M., and Ouhabaz, E.~M.}
\newblock Fractional anisotropic {C}alder\'{o}n problem on complete {R}iemannian manifolds.
\newblock {\em Commun. Contemp. Math. 26}, 9 (2024), Paper No. 2350057, 17.

\bibitem{Cov20}
{\sc Covi, G.}
\newblock Inverse problems for a fractional conductivity equation.
\newblock {\em Nonlinear Anal. 193\/} (2020), 111418, 18.

\bibitem{CGR22}
{\sc Covi, G., Garc\'ia-Ferrero, M.~A., and R\"uland, A.}
\newblock On the {C}alder\'on problem for nonlocal {S}chr\"odinger equations with homogeneous, directionally antilocal principal symbols.
\newblock {\em J. Differential Equations 341\/} (2022), 79--149.

\bibitem{CGRU23}
{\sc Covi, G., Ghosh, T., R\"{u}land, A., and Uhlmann, G.}
\newblock {A} reduction of the fractional {C}alder{\'{o}}n problem to the local {C}alder{\'{o}}n problem by means of the {C}affarelli-{S}ilvestre extension.
\newblock {\em arXiv:2305.04227\/} (2023).

\bibitem{Das25}
{\sc Das, S.}
\newblock {B}oundary control and {C}alder{\'{o}}n type inverse problems in non-local heat equation.
\newblock {\em arXiv:2504.20517\/} (2025).

\bibitem{DKSU09}
{\sc Dos Santos~Ferreira, D., Kenig, C.~E., Salo, M., and Uhlmann, G.}
\newblock Limiting {C}arleman weights and anisotropic inverse problems.
\newblock {\em Invent. Math. 178}, 1 (2009), 119--171.

\bibitem{DKLS16}
{\sc Dos Santos~Ferreira, D., Kurylev, Y., Lassas, M., and Salo, M.}
\newblock The {C}alder\'on problem in transversally anisotropic geometries.
\newblock {\em J. Eur. Math. Soc. (JEMS) 18}, 11 (2016), 2579--2626.

\bibitem{Feiz24}
{\sc Feizmohammadi, A.}
\newblock Fractional {C}alder\'{o}n problem on a closed {R}iemannian manifold.
\newblock {\em Trans. Amer. Math. Soc. 377}, 4 (2024), 2991--3013.

\bibitem{FGKRSU25}
{\sc Feizmohammadi, A., Ghosh, T., Krupchyk, K., R\"{u}land, A., Sj\"{o}strand, J., and Uhlmann, G.}
\newblock {F}ractional anisotropic {C}alder{\'{o}}n problem with external data.
\newblock {\em arXiv:2502.00710\/} (2025).

\bibitem{FGKU25}
{\sc Feizmohammadi, A., Ghosh, T., Krupchyk, K., and Uhlmann, G.}
\newblock {F}ractional anisotropic {C}alder{\'{o}}n problem on closed {R}iemannian manifolds.
\newblock {\em arXiv:2112.03480\/} (2021).

\bibitem{FKU24}
{\sc Feizmohammadi, A., Krupchyk, K., and Uhlmann, G.}
\newblock {C}alder{\'{o}}n problem for fractional {S}chr{\"{o}}dinger operators on closed {R}iemannian manifolds.
\newblock {\em arXiv:2407.16866\/} (2024).

\bibitem{FLW}
{\sc Frank, R.~L., Laptev, A., and Weidl, T.}
\newblock An improved one-dimensional {H}ardy inequality.
\newblock {\em J. Math. Sci. (N.Y.) 268}, 3, Problems in mathematical analysis. No. 118 (2022), 323--342.

\bibitem{Gho22}
{\sc Ghosh, T.}
\newblock A non-local inverse problem with boundary response.
\newblock {\em Rev. Mat. Iberoam. 38}, 6 (2022), 2011--2032.

\bibitem{GLX17}
{\sc Ghosh, T., Lin, Y.-H., and Xiao, J.}
\newblock The {C}alder\'on problem for variable coefficients nonlocal elliptic operators.
\newblock {\em Comm. Partial Differential Equations 42}, 12 (2017), 1923--1961.

\bibitem{GRSU20}
{\sc Ghosh, T., R\"uland, A., Salo, M., and Uhlmann, G.}
\newblock Uniqueness and reconstruction for the fractional {C}alder\'on problem with a single measurement.
\newblock {\em J. Funct. Anal. 279}, 1 (2020), 108505, 42.

\bibitem{GSU20}
{\sc Ghosh, T., Salo, M., and Uhlmann, G.}
\newblock The {C}alder\'{o}n problem for the fractional {S}chr\"{o}dinger equation.
\newblock {\em Anal. PDE 13}, 2 (2020), 455--475.

\bibitem{GU21}
{\sc Ghosh, T., and Uhlmann, G.}
\newblock {T}he {C}alder{\'{o}}n problem for nonlocal operators.
\newblock {\em arXiv:2110.09265\/} (2021).

\bibitem{GO08}
{\sc Gilboa, G., and Osher, S.}
\newblock Nonlocal operators with applications to image processing.
\newblock {\em Multiscale Model. Simul. 7}, 3 (2008), 1005--1028.

\bibitem{GS94}
{\sc Grigis, A., and Sj\"{o}strand, J.}
\newblock {\em Microlocal analysis for differential operators}, vol.~196 of {\em London Mathematical Society Lecture Note Series}.
\newblock Cambridge University Press, Cambridge, 1994.
\newblock An introduction.

\bibitem{grigor97}
{\sc Grigor’yan, A.}
\newblock Gaussian upper bounds for the heat kernel on arbitrary manifolds.
\newblock {\em J. Diff. Geom 45}, 1 (1997), 33--52.

\bibitem{GT11}
{\sc Guillarmou, C., and Tzou, L.}
\newblock Calder\'{o}n inverse problem with partial data on {R}iemann surfaces.
\newblock {\em Duke Math. J. 158}, 1 (2011), 83--120.

\bibitem{HLW24}
{\sc Harrach, B., Lin, Y., and {W}eth, T.}
\newblock {T}he {C}alder{\'{o}}n problem for the logarithmic ´ {S}chr{\"{o}}dinger equation.
\newblock {\em arXiv:2412.17775\/} (2024).

\bibitem{HL20}
{\sc Harrach, B., and Lin, Y.-H.}
\newblock Monotonicity-based inversion of the fractional {S}ch\"odinger equation {II}. {G}eneral potentials and stability.
\newblock {\em SIAM J. Math. Anal. 52}, 1 (2020), 402--436.

\bibitem{HLOS18}
{\sc Helin, T., Lassas, M., Oksanen, L., and Saksala, T.}
\newblock Correlation based passive imaging with a white noise source.
\newblock {\em J. Math. Pures Appl. (9) 116\/} (2018), 132--160.

\bibitem{KSU07}
{\sc Kenig, C.~E., Sj\"{o}strand, J., and Uhlmann, G.}
\newblock The {C}alder\'{o}n problem with partial data.
\newblock {\em Ann. of Math. (2) 165}, 2 (2007), 567--591.

\bibitem{KLW22}
{\sc Kow, P.-Z., Lin, Y.-H., and Wang, J.-N.}
\newblock The {C}alder\'{o}n problem for the fractional wave equation: uniqueness and optimal stability.
\newblock {\em SIAM J. Math. Anal. 54}, 3 (2022), 3379--3419.

\bibitem{LAWT21}
{\sc Laptev, A., and Weth, T.}
\newblock Spectral properties of the logarithmic {L}aplacian.
\newblock {\em Anal. Math. Phys. 11}, 3 (2021), Paper No. 133, 24.

\bibitem{LTU03}
{\sc Lassas, M., Taylor, M., and Uhlmann, G.}
\newblock The {D}irichlet-to-{N}eumann map for complete {R}iemannian manifolds with boundary.
\newblock {\em Comm. Anal. Geom. 11}, 2 (2003), 207--221.

\bibitem{LU01}
{\sc Lassas, M., and Uhlmann, G.}
\newblock On determining a {R}iemannian manifold from the {D}irichlet-to-{N}eumann map.
\newblock {\em Ann. Sci. \'{E}cole Norm. Sup. (4) 34}, 5 (2001), 771--787.

\bibitem{LAX02}
{\sc Lax, P.~D.}
\newblock {\em Functional analysis}.
\newblock Pure and Applied Mathematics (New York). Wiley-Interscience [John Wiley \& Sons], New York, 2002.

\bibitem{LU89}
{\sc Lee, J.~M., and Uhlmann, G.}
\newblock Determining anisotropic real-analytic conductivities by boundary measurements.
\newblock {\em Comm. Pure Appl. Math. 42}, 8 (1989), 1097--1112.

\bibitem{Li20}
{\sc Li, L.}
\newblock The {C}alder\'{o}n problem for the fractional magnetic operator.
\newblock {\em Inverse Problems 36}, 7 (2020), 075003, 14.

\bibitem{Nach96}
{\sc Nachman, A.~I.}
\newblock Global uniqueness for a two-dimensional inverse boundary value problem.
\newblock {\em Ann. of Math. (2) 143}, 1 (1996), 71--96.

\bibitem{PW87}
{\sc Paley, R. E. A.~C., and Wiener, N.}
\newblock {\em Fourier transforms in the complex domain}, vol.~19 of {\em American Mathematical Society Colloquium Publications}.
\newblock American Mathematical Society, Providence, RI, 1987.
\newblock Reprint of the 1934 original.

\bibitem{SP25}
{\sc {P}RAMANIK, S.}
\newblock {A}nisotropic {C}alder{\'{o}}n problem for a non-local second order elliptic operator.
\newblock {\em arXiv:2505.12255\/} (2025).

\bibitem{HU24}
{\sc Quan, H., and Uhlmann, G.}
\newblock The {C}alder\'{o}n problem for the fractional {D}irac operator.
\newblock {\em Math. Res. Lett. 31}, 1 (2024), 279--302.

\bibitem{Rudin64}
{\sc Rudin, W.}
\newblock {\em Principles of mathematical analysis}, second~ed.
\newblock McGraw-Hill Book Co., New York, 1964.

\bibitem{RS20}
{\sc R\"{u}land, A., and Salo, M.}
\newblock The fractional {C}alder\'{o}n problem: low regularity and stability.
\newblock {\em Nonlinear Anal. 193\/} (2020), 111529, 56.

\bibitem{RS2020}
{\sc R\"{u}land, A., and Salo, M.}
\newblock Quantitative approximation properties for the fractional heat equation.
\newblock {\em Math. Control Relat. Fields 10}, 1 (2020), 1--26.

\bibitem{Sal13}
{\sc Salo, M.}
\newblock The {C}alder\'on problem on {R}iemannian manifolds.
\newblock In {\em Inverse problems and applications: inside out. {II}}, vol.~60 of {\em Math. Sci. Res. Inst. Publ.} Cambridge Univ. Press, Cambridge, 2013, pp.~167--247.

\bibitem{SU87}
{\sc Sylvester, J., and Uhlmann, G.}
\newblock A global uniqueness theorem for an inverse boundary value problem.
\newblock {\em Ann. of Math. (2) 125}, 1 (1987), 153--169.

\bibitem{Tay23}
{\sc Taylor, M.~E.}
\newblock {\em Partial differential equations {I}. {B}asic theory}, vol.~115 of {\em Applied Mathematical Sciences}.
\newblock Springer, Cham, [2023] \copyright 2023.
\newblock Third edition [of 1395148].

\bibitem{Uhl14}
{\sc Uhlmann, G.}
\newblock Inverse problems: seeing the unseen.
\newblock {\em Bull. Math. Sci. 4}, 2 (2014), 209--279.

\bibitem{Zim23}
{\sc Zimmermann, P.}
\newblock Inverse problem for a nonlocal diffuse optical tomography equation.
\newblock {\em Inverse Problems 39}, 9 (2023), Paper No. 094001, 25.

\end{thebibliography}

\end{document}